\documentstyle[11pt,amssymb,amstex]{amsart}

\numberwithin{equation}{section}

\theoremstyle{plain}
\newtheorem{thm}{Theorem}[section]

\newtheorem{lem}[thm]{Lemma}
\newtheorem{prop}[thm]{Proposition}

\theoremstyle{definition}

\theoremstyle{remark}
\newtheorem{rem}{Remark}[section]

\numberwithin{equation}{section}

\def\bbf{{\mathbb F}}

\def\bn{{\mathbb N}}

\def\bq{{\mathbb Q}}
\def\br{{\mathbb R}}

\def\bz{{\mathbb Z}}

\begin{document}

\title[On cubic equations]
{On cubic equations over $P-$adic field}

\author{Farrukh Mukhamedov}
\address{Farrukh Mukhamedov\\
Department of Computational \& Theoretical Sciences \\
Faculty of Sciences, International Islamic University Malaysia\\
P.O. Box, 141, 25710, Kuantan\\
Pahang, Malaysia} \email{{\tt far75m@@yandex.ru}}

\author{Bakhrom Omirov}
\address{Bakhrom Omirov\\
Institute of Mathematics and Information Technologies, Tashkent, Uzbekistan}
\email{{\tt omirovb@@mail.ru}}

\author{Mansoor Saburov}
\address{Mansoor Saburov\\
Department of Computational \& Theoretical Sciences \\
Faculty of Science, International Islamic University Malaysia\\
P.O. Box, 141, 25710, Kuantan\\
Pahang, Malaysia} \email{{\tt msaburov@@gmail.com}}

\begin{abstract}
We provide a solvability criteria for a depressed cubic equation
in domains $\bz_p^{*},\bz_p,\bq_p$. We show that,
in principal, the Cardano method is not always applicable for such equations. Moreover, the numbers of solutions of the depressed cubic equation in domains $\bz_p^{*},\bz_p,\bq_p$ are provided. Since $\bbf_p\subset\bq_p,$ we generalize J.-P. Serre's \cite{JPSJ} and Z.H.Sun's \cite{ZHS1,ZHS3} results concerning  with depressed cubic equations over the finite field $\bbf_p$. Finally, all depressed cubic equations, for which the Cardano method could be applied, are described and the $p-$adic Cardano formula is provided for those cubic equations.


\end{abstract}

\maketitle

\section{Introduction}

In principal, according to Ostrowski's theorem (see \cite{FGou},\cite{NK},\cite{WS}), there are only two types of absolute values on the field of rational numbers $\bq$ : Archimedean or the real absolute value $|\cdot|_{\infty}$ and non-Archimedean or the $p-$adic absolute value $|\cdot|_p$ for some prime number $p$. These topological differences influence algebraic structures of the real and $p-$adic number fields $\br\equiv\overline{\bq}^{|\cdot|_{\infty}}$, $\bq_p\equiv\overline{\bq}^{|\cdot|_p}$.

Over the last century, $p-$adic numbers and $p-$adic analysis have come to
play a central role in modern number theory. This importance comes from
the fact that they afford a natural and powerful language for talking about
congruences between integers, and allow the use of methods borrowed from
analysis for studying such problems.

The fields of $p-$adic numbers were introduced by German
mathematician K. Hensel \cite{Hen}. The $p-$adic numbers were
motivated primarily by an attempt to bring the ideas and techniques
of the power series into number theory. Their canonical
representation is analogous to the expansion of analytic functions
into power series. This is one of the manifestations of the analogy
between algebraic numbers and algebraic functions.

For a fixed prime $p$, by $\bq_p$ it is denoted the field of
$p-$adic numbers, which is a completion of the rational numbers
$\bq$ with respect to the non-Archimedean norm $|\cdot|_p:\bq\to\br$
given by
\begin{eqnarray}
|x|_p=\left\{
\begin{array}{c}
  p^{-r} \ x\neq 0,\\
  0,\ \quad x=0,
\end{array}
\right.
\end{eqnarray}
here, $x=p^r\frac{m}{n}$ with $r,m\in\bz,$ $n\in\bn$,
$(m,p)=(n,p)=1$. A number $r$ is called \textit{a $p-$order} of $x$
and it is denoted by $ord_p(x)=r.$

Any $p-$adic number $x\in\bq_p$ can be uniquely represented in the
following canonical form
\begin{eqnarray*}
x=p^{ord_p(x)}\left(x_0+x_1\cdot p +x_2\cdot p^2+\cdots \right)
\end{eqnarray*}
where $x_0\in \{1,2,\cdots p-1\}$ and $x_i\in\{0,1,2,\cdots p-1\}$,
$i\geq 1,$ (see \cite{Bor Shaf}, \cite{NK})

More recently, numerous applications of $p-$adic numbers have shown
up in theoretical physics and quantum mechanics (see for example,
\cite{ArafDragFramVol}, \cite{BelGas}, \cite{FreWit}, \cite{Khren91,
Khren94}, \cite{Man}-\cite{MR2}, \cite{VladVolZel,Vol}).

The $p-$adic numbers are connected with solutions of Diophantine
equations modulo increasing powers of a prime number. The study of
Diophantine equations is finding solutions of polynomial equations
or systems of equations in integers, rational numbers, or sometimes
more general number rings. Such a topic is one of the oldest
branches of number theory, in fact of mathematics itself. The theory
of Diophantine equations in number rings was well developed in
\cite{Bor Shaf}, \cite{Con}.

One of the simplest Diophantine equation is the following equation
\begin{eqnarray}\label{eqxq=a}
x^q=a
\end{eqnarray}
over $\bq_p$, where $q\in\bn,$ $a\in\bq_p$. The solvability criterion for the equation \eqref{eqxq=a} from algebraic number theory point of view  was provided in \cite{SL}, \cite{JN}, \cite{JPS}. However, surprisingly, this criterion was not mentioned in the Bible books of the $p-$adic analysis (see \cite{FGou}, \cite{NK}, \cite{WS}) except $q=2$. From $p-$adic analysis point of view, the solvability criterion for the equation \eqref{eqxq=a} for any $q\in\bn$ was provided in \cite{AIRR1,AIRR2}, \cite{COR}, \cite{FMMS}.

All of us are aware that there is a criterion for the solvability of any quadratic equation $ax^2+bx+c=0$ in $\bq_p$,
where $a,b,c\in \bq_p$ with $a\neq 0$. It can be derived by the method of completing the square as follows. Without loss any generality, we may consider the following quadratic equation $x^2+qx+r=0,$ where $q=\frac{b}{a},$ $r=\frac{c}{a}.$ We then have that $(x+\frac{q}{2})^2=\left(\frac{q}{2}\right)^2-r.$ This equation has a solution in $\bq_p$ if and only if $\log\limits_p\left|\left(\frac{q}{2}\right)^2-r\right|_p$ is an even number and $\left(\left|\left(\frac{q}{2}\right)^2-r\right|_p\left(\left(\frac{q}{2}\right)^2-r\right)\right)^{\frac{p-1}{2}}\equiv 1 \ (mod \ p)$ for $p>2$ or $\left|\left(\frac{q}{2}\right)^2-r\right|_2\left(\left(\frac{q}{2}\right)^2-r\right)\equiv1 \ (mod \ 8)$ for $p=2.$ Thus, the solutions of the quadratic equation can be given in the form of $x_{\pm}=-\frac{q}{2}\pm\sqrt{\left(\frac{q}{2}\right)^2-r}.$

Unlike real numbers $\br,$ in general, the cubic equation $ax^3+bx^2+cx+d=0$ is not necessary to have a solution in $\bq_p,$ where $a,b,c,d\in \bq_p$ with $a\neq 0$. For example, the following simple cubic equation $x^3=p$ does not have any solution in $\bq_p.$ Therefore, it is a natural to find some criterion for solvability of the cubic equation in $\bq_p$. One of the general way to find solutions of the cubic equation in a local field is the Cardano method. However, by means of the Cardano method, we could not tell an existence of solutions of all cubic equations. Let us illustrate it in the following example.

We consider the following cubic equation $x^3-\frac{3}{p}x+\frac{p-3}{p}=0$ in $\bq_p$, where $p>3.$ It is clear that this cubic equation has a solution $x_{*}=-1.$ However, the Cardano method is not applicable for this equation. In fact, let us search for a solution in the form of $x=u+v.$ After elementary calculations, we obtain the following system of equations
$$
\left\{
\begin{array}{c}
uv=\frac{1}{p}\\
u^3+v^3=-\frac{p-3}{p}
\end{array}
\right.
$$
In order to solve this system, we should solve the following quadratic equation $z^2+\frac{p-3}{p}z+\frac{1}{p^3}=0$ in $\bq_p$, where $z=u^3$. However, $\log_p\left|\left(\frac{p-3}{2p}\right)^2-\frac{1}{p^3}\right|_p=3$ is odd. Thus, this quadratic equation $z^2+\frac{p-3}{p}z+\frac{1}{p^3}=0$  does not have solutions in $\bq_p$. Therefore, in general, by means of Cardano method we could not detect solutions of all cubic equations.

To the best of our knowledge, we could not find the solvability criterion in an explicit form for the cubic equation \eqref{cubiceqn} in the Bible books of $p-$adic analysis and algebraic number theory (see  \cite{Apos}, \cite{FGou}, \cite{NK}, \cite{SL}, \cite{JN}, \cite{WS}, \cite {JPS}). Some sufficient conditions for the solvability of the cubic equation \eqref{cubiceqn} in $\bq_p$ were provided in \cite{AIRR1,AIRR2}, \cite{JN}.

In this paper, we provide the criterion for solvability of any cubic equation $ax^3+bx^2+cx+d=0$ in $\bq_p$ where $a,b,c,d\in \bq_p$ with $a\neq 0$. Dividing by $a$ and substituting $x$ by $x-\frac{b}{3a}$ we can get the so-called \emph{depressed cubic} equation
\begin{eqnarray}\label{cubiceqn}
x^3+ax=b
\end{eqnarray}
Concerning with classification problems of finite dimensional Leibniz algebras (see \cite{KhudKurb}), in this paper, we are going to provide the criterion for the solvability of the depressed cubic equation \eqref{cubiceqn}, where $a,b\in \bq_p,$ in domains $\bz_p^{*}$, $\bz_p$, $\bq_p$. In general case, one can easily derive the criterion from the depressed cubic equation.

It is worth mentioning that there are some cubic equations which do not have any solutions in $\bz_p^{*}$ (in $\bz_p$) but have solutions in $\bz_p$ (in $\bq_p$).

Let us consider the following equation $x^3+p^2x=2p^3$ in $\bq_p.$ This equation does not have any solutions in $\bz_p^{*}.$ In fact, for any $x\in\bz_p^{*}$ one has that $\left|x^3+p^2x\right|_p=1$. On the other hand, $|2p^3|_p<1$. This contradictions shows that the cubic equation $x^3+p^2x=2p^3$ does not have any solution in $\bz_p^{*}.$ However, it has a solution $x_{*}=p$ which belongs to $\bz_p.$

Let us consider the following equation $x^3+p^2x=\frac{1+p^4}{p^3}$ in $\bq_p.$ If $x\in\bz_p$ then $|x^3+p^2x|_p\leq 1.$ However, $\left|\frac{1+p^4}{p^3}\right|_p>1$. This means that the equation $x^3+p^2x=\frac{1+p^4}{p^3}$ does not have any solution in $\bz_p.$ On the other hand, this equation has a solution $x_{*}=\frac{1}{p}$ which belongs to $\bq_p.$

Therefore, finding the criterion for the solvability of the depressed cubic equation \eqref{cubiceqn}, where $a,b\in \bq_p,$ in domains $\bz_p^{*}$, $\bz_p$, $\bq_p$ is of independent interest.

Let us consider the depressed cubic equation \eqref{cubiceqn} in the finite field $\bbf_p= \bz\diagup p\bz$
\begin{eqnarray}\label{cubcongr}
x^3+ax=b,
\end{eqnarray}
where $a,b,x\in \bbf_p.$ By ${\mathbf{N}}_{\bbf_p}(x^3+ax-b)$, we denote the number of solutions of the depressed cubic equation \eqref{cubcongr} in $\bbf_p$, where $a,b\in \bbf_p$.

If $a=\bar{0}$ or $b=\bar{0}$ then ${\mathbf{N}}_{\bbf_p}(x^3+ax-b)$ can be easily study. Therefore, we shall assume that $a\neq\bar{0}$ and $b\neq\bar{0}$.

The study of ${\mathbf{N}}_{\bbf_p}(x^3+ax-b)$ goes back to several centuries with contributions from  Euler, Gauss, Jacobi, Cauchy. Surprisingly, up to now, there are still many papers which were devoted to study ${\mathbf{N}}_{\bbf_p}(x^3+ax-b)$ (see papers \cite{JPSJ}, \cite{ZHS1,ZHS3}, and references therein)

We denote the number of solutions of the depressed cubic equation \eqref{cubiceqn}, where  $a,b\in \bq_p,$  in domains $\bz_p^{*}$, $\bz_p$, $\bq_p$ respectively by ${\mathbf{N}}_{\bz_p^{*}}(x^3+ax-b)$, ${\mathbf{N}}_{\bz_p}(x^3+ax-b)$, ${\mathbf{N}}_{\bq_p}(x^3+ax-b)$. In this paper we shall give the description of all these sets whenever $a,b\in \bq_p$. Since $\bbf_p\subset\bq_p$, our
results generalize all results  in \cite{JPSJ}, \cite{ZHS1,ZHS3} which were concerning with equation \eqref{cubcongr} in $\bbf_p$.

The last but not least, as we already mentioned that the Cardano method is not applicable for all depressed cubic equations \eqref{cubiceqn}, where $a,b\in \bq_p$. Therefore, in the last section we are going to describe all depressed cubic equations for which the Cardano method can be applied and we shall provide $p-$adic Cardano formula for those equations.

\section{Preliminaries}
In this section we shall recall some necessary results form number theory.

We respectively denote the set of all {\it $p-$adic integers} and {\it units} of
$\bq_p$ by
$$\bz_p=\{x\in\bq_{p}: |x|_p\leq1\}, \quad \bz_p^{*}=\{x\in\bq_{p}: |x|_p=1\}.$$

Any $p-$adic unit $x\in\bz_p^{*}$ has the following unique canonical form
$$x=x_0+x_1\cdot p+x_2\cdot p^2+\cdots$$
where $x_0\in \{1,2,\cdots p-1\}$ and $x_i\in\{0,1,2,\cdots p-1\}$, $i\in\bn.$

Any nonzero $p-$adic number $x\in\bq_p$ has a unique representation of the form $x =\cfrac{x^{*}}{|x|_p}$, where $x^{*}\in\bz_p^{*}$.

\begin{lem}[Hensel's Lemma, \cite{Bor Shaf}]\label{Hensel}
Let $f(x)$ be polynomial whose the coefficients are $p-$adic
integers. Let $\theta$ be a $p-$adic integer such that for some
$i\geq 0$ we have
$$
f(\theta)\equiv 0 \ (mod \ p^{2i+1}),
$$
$$
f'(\theta)\equiv 0 \ (mod \ p^{i}), \quad f'(\theta)\not\equiv 0 \ (mod \ p^{i+1}).
$$
Then $f(x)$ has a unique $p-$adic integer root $x_0$ which satisfies $x_0\equiv \theta\ (mod \ p^{i+1}).$
\end{lem}

Let $p$ be a prime number, $q\in\bn$, $a\in\bbf_p$ with $a\neq \bar{0}.$  The number $a$ is called \textit{a $q$-th power
residue modulo $p$} if the the following equation
\begin{eqnarray}\label{kthresidue}
x^q=a
\end{eqnarray} has a solution in $\bbf_p$.

\begin{prop}[\cite{Ros}]\label{aisresidueofp}
Let $p$ be an odd prime number, $q\in\bn$, $d=(q,p-1),$ and $a\in\bbf_p$ with $a\neq\bar{0}.$ Then the following statements hold true:
\begin{itemize}
  \item [(i)] $a$ is the $q$-th power residue modulo $p$ if and only if one has
$a^{\frac{p-1}{d}}=\bar{1};$
  \item [(ii)] If $a^{\frac{p-1}{d}}=\bar{1}$ then the equation \eqref{kthresidue} has $d$ number of solutions in $\bbf_p$.
\end{itemize}
\end{prop}

The solvability criterion of the following equation in $\bq_p$
\begin{eqnarray}\label{x^q=a}
x^q=a,
\end{eqnarray}
where $q\in\bn$, $a\in \bq_p$ with $a\neq 0$, was provided in  \cite{AIRR1,AIRR2}, \cite{COR}, \cite{SL}, \cite{FMMS}, \cite{JN}, \cite{JPS}.
\begin{prop}\label{Criterionforp}
Let $p$ be an odd prime number, $q\in\bn,$ $a\in \bq_p,$ $a=\frac{a^{*}}{|a|_p}$ and $a^{*}\in\bz_p^{*}$ with $a^{*}=a_0+a_1\cdot p+a_2\cdot p^2+\cdots$. Then the following statements hold true:
\begin{itemize}
  \item [(i)] If $(q,p)=1$ then the equation \eqref{x^q=a} has a solution in $\bq_p$ if and only if $a_0^{\frac{p-1}{(q,p-1)}}\equiv 1 \ ( mod \ p)$ and $\log_p{|a|_p}$ is divisible by $q$
  \item [(ii)] If $(q,p)=1$, $a_0^{\frac{p-1}{(q,p-1)}}\equiv 1 \ ( mod \ p)$ and $\log_p{|a|_p}$ is divisible by $q$ then the equation \eqref{x^q=a} has $(q,p-1)$ number of solutions in $\bq_p$.
  \item [(iii)] If $(q,p)=1$ then the equation $y^q\equiv a^{*} \ (mod \ p^{k})$ has a solution in $\bz_p^{*}$ for some $k\in\bn$ if and only if $a_0^{\frac{p-1}{(q,p-1)}}\equiv 1 \ ( mod \ p).$
  \item [(iv)] If $q=m\cdot p^s$ with $(m,p)=1$, $s\geq 1$ then the equation \eqref{x^q=a}
  has a solution in $\bq_p$ if and only if $a_0^{\frac{p-1}{(m,p-1)}}\equiv 1 \ ( mod \ p)$, $a^{p^s}_0\equiv a \ (mod \ p^{s+1})$ and $\log_p{|a|_p}$ is divisible by $q$.
\end{itemize}
\end{prop}

Let us consider the following depressed cubic equation in the field $\bbf_p$
\begin{eqnarray}\label{cubiccong}
x^3+\bar{a}x=\bar{b},
\end{eqnarray}
where $\bar{a},\bar{b}\in \bbf_p.$  We assume that $\bar{a}\neq\bar{0}$ and $\bar{b}\neq\bar{0}$. The number of solutions ${\mathbf{N}}_{\bbf_p}(x^3+\bar{a}x-\bar{b})$ of this equation was described in \cite{ZHS2}.

\begin{prop}[\cite{ZHS2}]\label{CubicinF_p}
Let $p>3$ be a prime number and $\bar{a},\bar{b}\in\bbf_p$ with $\bar{a}\bar{b}\neq\bar{0}$. Let $\overline{D}=-4\bar{a}^3-27\bar{b}^2$ and $u_{n+3}=\bar{b}u_n-\bar{a}u_{n+1}$ for $n\in\bn$ with $u_1=\bar{0},$ $u_2=-\bar{a},$ $u_3=\bar{b}.$ Then the following holds true:
$$
{\mathbf{N}}_{\bbf_p}(x^3+\bar{a}x-\bar{b})=\left\{
\begin{array}{l}
3 \ \ \  if \ \ \ \overline{D}u_{p-2}^2=\bar{0} \\
0 \ \ \ if \ \ \ \overline{D}u_{p-2}^2=9\bar{a}^2 \\
1 \ \ \ if \ \ \ \overline{D}u_{p-2}^2\neq \bar{0}, 9\bar{a}^2
\end{array}
\right.
$$
\end{prop}
It is worth mentioning that ${\mathbf{N}}_{\bbf_p}(x^3+\bar{a}x-\bar{b})$ was studied in \cite{JPSJ} for the case $\bar{a}=-\bar{1}$ and $\bar{b}=\bar{1}$.

The proofs of the following statements are straightforward.
\begin{prop}\label{numberofcongequation}
Let $p>3$ be a prime number and $\bar{a},\bar{b}\in\bbf_p$, $\bar{a}\bar{b}\neq\bar{0}$. Let $\overline{D}=-4\bar{a}^3-27\bar{b}^2$ and $u_{n+3}=\bar{b}u_n-\bar{a}u_{n+1}$ with $\overline{D}u_{p-2}^2\neq9\bar{a}^2,$
$u_1=\bar{0},$ $u_2=-\bar{a},$ $u_3=\bar{b}.$
\begin{itemize}
  \item [I] Let $\overline{D}u_{p-2}^2=\bar{0}.$ Then the following statements holds true:
  \begin{itemize}
  \item [I.1] The equation \eqref{cubiccong} has 3 distinct solutions in $\bbf_p$ if and only if $\overline{D}\neq\bar{0}.$ Moreover, one has that $3\bar{x}^2+\bar{a}\neq0$ for any root $\bar{x};$
  \item [I.2] The equation \eqref{cubiccong} has 2 distinct solutions in $\bbf_p$ while one of them of multiplicity 2 if and only if $\bar{D}=\bar{0}.$ If $\bar{x}_1$, $\bar{x}_2$ are 2 distinct solutions while $\bar{x}_1$ is a multiple solution then $\bar{x}_1=\frac{3\bar{b}}{2\bar{a}},$ $\bar{x}_2=-\frac{3\bar{b}}{\bar{a}},$ and $3\bar{x}_2^2+\bar{a}\neq\bar{0};$
  \item [I.3] The equation \eqref{cubiccong} does not have any solution of multiplicity 3.
\end{itemize}

  \item [II] Let $\overline{D}u_{p-2}^2\neq\bar{0}, 9\bar{a}^2.$ If $\bar{x}$ is a solution of the equation \eqref{cubiccong} then $3\bar{x}^2+\bar{a}\neq\bar{0}.$
\end{itemize}
\end{prop}

\begin{rem} Due to Proposition \ref{numberofcongequation}, one may conclude that under the assumption of Proposition \ref{numberofcongequation}, there always exists at least one solution $\bar{x}$ of the equation \eqref{cubiccong} such that $3\bar{x}^2+\bar{a}\neq\bar{0}.$
\end{rem}

\section{The solvability criterion in domains $\bz_p^{*},$ $\bz_p,$ $\bq_p$ with $p>3$}

In this section we provide the solvability criterion for the depressed cubic equation \eqref{cubiceqn} in domains $\bz_p^{*},$ $\bz_p,$ $\bq_p$ with $p>3$, where $a,b\in\bq_p$ with $ab\neq 0$. The solvability criteria of the equation \eqref{cubiceqn} for the case $ab=0$ is given in \cite{COR}, \cite{FMMS} (see Proposition \ref{Criterionforp}).

We need the following auxiliary result.

\begin{prop}\label{NecessaryconditionforZ_p^*}
Let $p$ be any prime. Suppose that the equation \eqref{cubiceqn} has a solution in $\bz_p^{*},$  where $a,b\in\bq_p$ with $ab\neq 0.$ Then one of the following conditions holds true:
\begin{itemize}
  \item [(i)] $|a|_p<|b|_p=1;$
  \item [(ii)] $|b|_p<|a|_p=1;$
  \item [(iii)] $|a|_p=|b|_p\geq 1.$
\end{itemize}
\end{prop}
\begin{pf}
Let $p$ be any prime. We suppose that the equation \eqref{cubiceqn} has a solution in $\bz_p^{*}.$ Since $ab\neq 0,$
one can get that
\begin{eqnarray*}
|b|_p=|x^3+ax|_p\leq\max\{1,|a|_p\}, \\
|a|_p=|ax|_p=|b-x^3|_p\leq\max\{1,|b|_p\},\\
1=|x^3|_p=|b-ax|_p\leq\max\{|a|_p,|b|_p\}.
\end{eqnarray*}
Thus, if $|a|_p\neq |b|_p$ then $\max\{|a|_p,|b|_p\}=1$ and if $|a|_p=|b|_p$ then $|a|_p=|b|_p\geq1.$ This yields the claim.
\end{pf}

We are going to state the solvability criterion for the depressed cubic equation \eqref{cubiceqn} in domains $\bz_p^{*},$ $\bz_p,$ $\bq_p$ for $p>3.$

Let $a,b\in\bq_p$ be two nonzero $p-$adic numbers with $a=\frac{a^{*}}{|a|_p},$ $b=\frac{b^{*}}{|b|_p}$ where $a^{*},b^{*}\in\bz_p^{*}$ with $a^{*}=a_0+a_1\cdot p+a_2\cdot p^2+\cdots$ and $b^{*}=b_0+b_1\cdot p+b_2\cdot p^2+\cdots.$

We set $D_0=-4a_0^3-27b_0^2$ and $u_{n+3}=b_0u_n-a_0u_{n+1}$ with $u_1=0,$ $u_2=-a_0,$ and $u_3=b_0$ for $n=\overline{1,p-3}$

\begin{thm}\label{criterionforp>3}
Let $p>3$ be a prime. Then the following statements hold true:

\begin{itemize}
  \item [I]
  The equation \eqref{cubiceqn} has a solution in $\bz_p^{*}$ if and only if one of the following conditions holds true:
  \begin{itemize}
    \item [I.1] $|a|_p<|b|_p=1$ and $b_0^{\frac{p-1}{(3,p-1)}}\equiv 1 \ (mod \ p);$
    \item [I.2] $|b|_p<|a|_p=1$ and $(-a_0)^{\frac{p-1}{2}}\equiv 1 \ (mod \ p);$
    \item [I.3] $|a|_p=|b|_p=1$ and $D_0u_{p-2}^2\not\equiv 9a_0^{2} \ (mod \ p);$
    \item [I.4] $|a|_p=|b|_p>1.$
  \end{itemize}
  \item [II]
  The equation \eqref{cubiceqn} has a solution in $\bz_p$ if and only if one of the following conditions holds true:
  \begin{itemize}
    \item [II.1] $|a|_p^3<|b|_p^2\leq 1,$ \ $3\mid\log_p|b|_p,$ and $b_0^{\frac{p-1}{(3,p-1)}}\equiv 1 \ (mod \ p);$
    \item [II.2] $|a|_p^3=|b|_p^2\leq 1$ and $D_0u_{p-2}^2\not\equiv 9a_0^{2} \ (mod \ p);$
    \item [II.3] $|a|_p^3>|b|_p^2$ and $|a|_p\geq |b|_p.$
  \end{itemize}
  \item [III]
  The equation \eqref{cubiceqn} has a solution in $\bq_p$ if and only if one of the following conditions holds true:
  \begin{itemize}
    \item [III.1] $|a|_p^3<|b|_p^2,$ \ $3\mid\log_p|b|_p,$ and $b_0^{\frac{p-1}{(3,p-1)}}\equiv 1 \ (mod \ p);$
    \item [III.2] $|a|_p^3=|b|_p^2$ and $D_0u_{p-2}^2\not\equiv 9a_0^{2} \ (mod \ p);$
    \item [III.3] $|a|_p^3>|b|_p^2.$
  \end{itemize}
\end{itemize}
\end{thm}
\begin{pf}
I. Let $p>3$ be a prime number and $a,b\in\bq_p$ be two nonzero $p-$adic numbers. Due to Proposition \ref{NecessaryconditionforZ_p^*}, if  the equation \eqref{cubiceqn} has a solution in $\bz_p^{*}$ then one of the following assertions holds true: (i) $|a|_p<1,$  $|b|_p=1;$ or
(ii) $|a|_p=1,$ $|b|_p<1;$  or (iii) $|a|_p=|b|_p\geq 1.$ Let us study every case.

I.1.  Now, we want to show that if  the equation \eqref{cubiceqn} has a solution in $\bz_p^{*}$ then under the condition $|a|_p<1,$  $|b|_p=1$ one has that $b_0^{\frac{p-1}{(3,p-1)}}\equiv 1 \ (mod \ p).$ In fact, let $x\in\bz_p^{*}$ be a solution of the equation \eqref{cubiceqn}. Since $|a|_p<1$ and $|b|_p=1,$ we get that
$$
x_0^3+ax_0\equiv x_0^3 \equiv b_0 \ (mod \ p).
$$
Due to Proposition \ref{aisresidueofp}, this yields that $b_0^{\frac{p-1}{(3,p-1)}}\equiv 1 \ (mod \ p).$

Now, we are going to prove that if  $|a|_p<1,$  $|b|_p=1,$ and $b_0^{\frac{p-1}{(3,p-1)}}\equiv 1 \ (mod \ p)$ then the equation \eqref{cubiceqn} has a solution in $\bz_p^{*}$. In fact, since $b_0^{\frac{p-1}{(3,p-1)}}\equiv 1 \ (mod \ p)$ there is $x_0\in\bz$ such that $x_0^3\equiv b_0 \ (mod \ p)$ and $(x_0,p)=1$. Let us consider the following polynomial $f(x)=x^3+ax-b.$ Since $|a|_p<1,$ it is clear that
$$
f(x_0)=x_0^3+ax_0-b\equiv x_0^3-b_0\equiv 0 \ (mod \ p),\quad
f'(x_0)=3x_0^2+a\not\equiv 0 \ (mod \  p).
$$
Due to Hensel's Lemma \ref{Hensel} the equation \eqref{cubiceqn} has a solution $x\in\bz_p$ with $|x-x_0|_p\leq \frac{1}{p}$ where $|x_0|_p=1.$ This means that $|x|_p=1$, i.e., $x\in\bz_p^{*}.$

I.2. Let us show that if the equation \eqref{cubiceqn} has a solution in $\bz_p^{*}$ then under the condition $|a|_p=1,$ $|b|_p<1$ one has that $(-a_0)^{\frac{p-1}{2}}\equiv 1 \ (mod \ p).$ In fact, suppose $x\in\bz_p^{*}$ is a solution of the equation \eqref{cubiceqn}. Since $|a|_p=1$ and $|b|_p<1$  we get that
$$
x_0^3+a_0x_0\equiv x_0(x_0^2+a_0)\equiv b\equiv 0 \ (mod \ p).
$$
This means that $x_0^2+a_0\equiv 0 \ (mod \ p)$ and due to to Proposition \ref{aisresidueofp} we have that $(-a_0)^{\frac{p-1}{2}}\equiv 1 \ (mod \ p).$

Now, we want to prove that if  $|a|_p=1,$  $|b|_p<1,$ and $(-a_0)^{\frac{p-1}{2}}\equiv 1 \ (mod \ p)$ then the equation \eqref{cubiceqn} has a solution in $\bz_p^{*}$. In fact, since $(-a_0)^{\frac{p-1}{2}}\equiv 1 \ (mod \ p)$ there is $x_0\in\bz$ such that $x_0^2+a_0\equiv 0 \ (mod \ p)$ and $(x_0,p)=1$. Let us consider the following polynomial $f(x)=x^3+ax-b.$ Since $|b|_p<1,$ it is clear that
\begin{eqnarray*}
f(x_0)&=&x_0^3+ax_0-b\equiv x_0(x_0^2+a_0)\equiv 0 \ (mod \ p),\\
f'(x_0)&=&3x_0^2+a\equiv3(x_0^2+a_0)-2a_0\not\equiv 0 \ (mod \  p).
\end{eqnarray*}
Due to Hensel's Lemma \ref{Hensel} the equation \eqref{cubiceqn} has a solution $x\in\bz_p$ with $|x-x_0|_p\leq \frac{1}{p}$ where $|x_0|_p=1.$ This yields that $|x|_p=1$, i.e., $x\in\bz_p^{*}.$

I.3. We are going to show that if the equation \eqref{cubiceqn} has a solution in $\bz_p^{*}$ then under the condition $|a|_p=1,$ $|b|_p=1$ we have that $D_0u_{p-2}^2\not\equiv 9a_0^{2} \ (mod \ p).$ In fact, let $x\in\bz_p^{*}$ be a solution of the equation \eqref{cubiceqn}. Since $|a|_p=1$ and $|b|_p=1$  we get that
$$
x_0^3+a_0x_0\equiv b_0 \ (mod \ p).
$$
This means that the depressed cubic equation $x^3+a_0x=b_0$ has at least one solution in the finite field $\bbf_p$. Due to to Proposition \ref{CubicinF_p} we have that $D_0u_{p-2}^2\not\equiv 9a_0^{2} \ (mod \ p).$

Now, we are going to prove that if $|a|_p=1,$  $|b|_p=1,$ and $D_0u_{p-2}^2\not\equiv 9a_0^{2} \ (mod \ p)$ then the equation \eqref{cubiceqn} has a solution in $\bz_p^{*}$. Since $D_0u_{p-2}^2\not\equiv 9a_0^{2} \ (mod \ p)$,  the equation $x^3+a_0x\equiv b_0 \ (mod \ p)$ has at least one solution. Due to Proposition \ref{numberofcongequation}, among all solutions of the equation $x^3+a_0x\equiv b_0 \ (mod \ p)$, there always exists at least one solution $x_0$ such that $3x_0^2+a_0\not\equiv 0 \ (mod \ p)$ and $(x_0,p)=1.$

Again, if we apply Hensel's Lemma \ref{Hensel} to the polynomial $f(x)=x^3+ax-b$ at the point $x=x_0$ then we get that the equation \eqref{cubiceqn} has a solution $x\in\bz_p$ with $|x-x_0|_p\leq \frac{1}{p}$ where $|x_0|_p=1.$ This yields that $|x|_p=1$, i.e., $x\in\bz_p^{*}.$

I.4. Now, we want to prove that if  $|a|_p=|b|_p>1$ then the equation \eqref{cubiceqn} has a solution in $\bz_p^{*}$. Let $|a|_p=|b|_p=p^m,$ where $m\geq 1.$ Then $a=p^{-m}a^{*},$ $b=p^{-m}b^{*}$ with $a^{*},b^{*}\in\bz_p^{*}.$ It is clear that the depressed equation $x^3+p^{-m}a^{*}x=p^{-m}b^{*}$ has a solution in $\bz_p^{*}$ if and only if the equation $p^{m}x^3+a^{*}x=b^{*}$ has a solution in $\bz_p^{*}$, where $a^{*},b^{*}\in\bz_p^{*}.$ To this end, we consider the following polynomial $f(x)=p^{m}x^3+a^{*}x-b^{*}.$ If $a_0,b_0$ are the first digits of $a^{*},b^{*}\in\bz_p^{*}$ then there is $x_0\in\bz$ such that $a_0x_0\equiv b_0 \ (mod \ p)$ and $(x_0,p)=1.$ Then, we obtain that
$$
f(x_0)\equiv p^{m}x_0^3+a_0x_0-b_0\equiv 0 \ (mod \ p), \quad
f'(x_0)\equiv 3p^{m}x_0^2+a_0\not\equiv 0\ (mod \ p).
$$
Again, due to Hensel's Lemma \ref{Hensel}, the equation $p^{m}x^3+a^{*}x=b^{*}$ has a solution $x\in\bz_p^{*}$ with $|x-x_0|_p\leq \frac{1}{p}$ where $|x_0|_p=1.$ This yields that $x\in\bz_p^{*}.$

II. Let $p>3$ be a prime number and $a,b\in\bq_p$ be two nonzero $p-$adic numbers. Now, we want to provide the solvability criterion for the depressed equation \eqref{cubiceqn} in the domain $\bz_p$.

We know that any $p-$adic integer $x$ has the following unique form $x=p^{k}x^{*},$ where $x^{*}\in\bz_p^{*}$ and $k\in\bn\cup\{0\}.$ Then the depressed equation \eqref{cubiceqn} has a solution in $\bz_p$ if and only if there is a nonnegative integer $k$ such that the following equation
\begin{eqnarray}\label{x^3AxB}
x^3+Ax=B
\end{eqnarray}
has a solution in $\bz_p^{*}$ where $A=p^{-2k}a$, $B=p^{-3k}b$.

Due to the case I, in general, the equation \eqref{x^3AxB} has a solution in $\bz_p^{*}$ if and only if one of the following conditions holds true:
\begin{itemize}
    \item [(i)] $|A|_p<|B|_p=1,$ and $b_0^{\frac{p-1}{(3,p-1)}}\equiv 1 \ (mod \ p);$
    \item [(ii)]   $|B|_p<|A|_p=1,$ and $(-a_0)^{\frac{p-1}{2}}\equiv 1 \ (mod \ p);$
    \item [(iii)] $|A|_p=|B|_p=1,$ and $D_0u_{p-2}^2\not\equiv 9a_0^{2} \ (mod \ p);$
    \item [(iv)] $|A|_p=|B|_p>1.$
\end{itemize}
where $A=p^{-2k}\frac{a^{*}}{|a|_p}$, $B=p^{-3k}\frac{b^{*}}{|b|_p}$, and $a^{*},b^{*}\in\bz_p^{*}$ with
$$a^{*}=a_0+a_1\cdot p+a_2\cdot p^2+\cdots; \quad b^{*}=b_0+b_1\cdot p+b_2\cdot p^2+\cdots.$$

It is clear that $|A|_p=p^{2k}|a|_p,$ $|B|_p=p^{3k}|b|_p.$ Now, in every case (i)-(iv), we want to describe all  $p-$adic numbers $a,b\in \bq_p$ for which the equation \eqref{x^3AxB} has a solution in $\bz_p^{*}$ for some nonnegative integer $k$.

II.1. Suppose that $|A|_p<|B|_p=1$ and $b_0^{\frac{p-1}{(3,p-1)}}\equiv 1 \ (mod \ p).$ Since $|B|_p=p^{3k}|b|_p,$ we get from $|B|_p=1$ that $3k=-\log_p|b|_p.$ The last equation has a nonnegative integer solution w.r.t $k$ if and only if $\log_p|b|_p$ is divisible by 3 and $|b|_p\leq 1$. Therefore, if $k=-\frac{\log_p|b|_p}{3}$ then it follows from $|A|_p=p^{2k}|a|_p<1$ that $|a|_p^3<|b|_p^2.$ Consequently, if $|a|_p^3<|b|_p^2\leq 1,$ $b_0^{\frac{p-1}{(3,p-1)}}\equiv 1 \ (mod \ p)$ and $3|\log_p|b|_p$ then the equation \eqref{x^3AxB} has a solution in $\bz_p^{*}$ for $k=-\frac{\log_p|b|_p}{3}$.

II.2. Assume that $|A|_p=|B|_p=1$ and $D_0u_{p-2}^2\not\equiv 9a_0^{2} \ (mod \ p).$ Since $|A|_p=p^{2k}|a|_p,$ $|B|_p=p^{3k}|b|_p$ we have that $2k=-\log_p|a|_p$ and $3k=-\log_p|b|_p.$ The last two equations have a nonnegative integer solution w.r.t $k$ if and only if $\log_p|a|_p$ is divisible by 2, $\log_p|b|_p$ is divisible by 3, and $|a|_p\leq 1$, $|b|_p\leq 1.$ Therefore, if $k=-\frac{\log_p|a|_p}{2}=-\frac{\log_p|b|_p}{3}$ then we have that $|a|_p^3=|b|_p^2.$ Consequently, if $|a|_p^3=|b|_p^2\leq 1$ and $D_0u_{p-2}^2\not\equiv 9a_0^{2} \ (mod \ p)$ then the equation \eqref{x^3AxB} has a solution in $\bz_p^{*}$ for $k=-\frac{\log_p|a|_p}{2}=-\frac{\log_p|b|_p}{3}$.

II.3. We are going to consider the cases $|B|_p<|A|_p=1,$ $(-a_0)^{\frac{p-1}{2}}\equiv 1 \ (mod \ p)$ and $|A|_p=|B|_p>1.$ Let us study every case.

Let $|B|_p<|A|_p=1$ and $(-a_0)^{\frac{p-1}{2}}\equiv 1 \ (mod \ p).$ Since $|A|_p=p^{2k}|a|_p$ we get from $|A|_p=1$ that $2k=-\log_p|a|_p.$ The last equation has a nonnegative integer solution w.r.t $k$ if and only if $\log_p|a|_p$ is even and $|a|_p\leq 1$. Therefore, if $k=-\frac{\log_p|a|_p}{2}$ then it follows from $|B|_p=p^{3k}|b|_p<1$ that $|a|_p^3>|b|_p^2.$ Consequently, if $|b|_p^2<|a|_p^3\leq1,$ $\log_p|a|_p$ is even, and $(-a_0)^{\frac{p-1}{2}}\equiv 1 \ (mod \ p)$ then the equation \eqref{x^3AxB} has a solution in $\bz_p^{*}$ for $k=-\frac{\log_p|a|_p}{2}$.

Let $|A|_p=|B|_p>1.$ Since $|A|_p=p^{2k}|a|_p,$ $|B|_p=p^{3k}|b|_p$ we obtain from $|A|_p=|B|_p$ that $k=\log_p|a|_p-\log_p|b|_p.$ Hence, $k$ is a nonnegative integer if and only if $|a|_p\geq|b|_p.$ Therefore, if $k=\log_p|a|_p-\log_p|b|_p$ then it follows form $|A|_p=|B|_p>1$ that $|a|_p^3>|b|_p^2.$ Consequently, if $|a|_p^3>|b|_p^2$ and $|a|_p\geq|b|_p$ then the equation \eqref{x^3AxB} has a solution in $\bz_p^{*}$ for $k=\log_p|a|_p-\log_p|b|_p$.

It is worth mentioning that if $p-$adic numbers $a,b\in\bq_p$ satisfy the conditions  $|b|_p^2<|a|_p^3\leq1,$ $\log_p|a|_p$ is even, and $(-a_0)^{\frac{p-1}{2}}\equiv 1 \ (mod \ p)$ then they satisfy the conditions $|a|_p^3>|b|_p^2$ and $|a|_p\geq|b|_p$ as well. Consequently, regardless of whether $p-$adic numbers $a,b\in\bq_p$ satisfy the conditions $|b|_p^2<|a|_p^3\leq1,$ $\log_p|a|_p$ is even, and $(-a_0)^{\frac{p-1}{2}}\equiv 1 \ (mod \ p)$, the equation \eqref{x^3AxB} has a solution in $\bz_p^{*}$ if $|a|_p^3>|b|_p^2$ and $|a|_p\geq|b|_p$. Moreover, if $|b|_p^2<|a|_p^3\leq1,$ $\log_p|a|_p$ is even, and $(-a_0)^{\frac{p-1}{2}}\equiv 1 \ (mod \ p)$ then the equation \eqref{x^3AxB} has at least two distinct solutions in $\bz_p^{*}$ for two distinct nonnegative integers $k=-\frac{\log_p|a|_p}{2}$ and $k=\log_p|a|_p-\log_p|b|_p,$ otherwise the equation \eqref{x^3AxB} has at least one solution in $\bz_p^{*}$ for $k=\log_p|a|_p-\log_p|b|_p$.

III. Let $p>3$ be a prime number and $a,b\in\bq_p$ be two nonzero $p-$adic numbers. In order to provide the solvability criterion for the depressed equation \eqref{cubiceqn} in the domain $\bq_p$, we apply the same method that one used in the case II, but with the assumption that $k$ is any integer number. This completes the proof.
\end{pf}

\section{The number of solutions in domains $\bz_p^{*},$ $\bz_p,$ $\bq_p,$ with $p>3$}

In this section, we are aiming to describe the numbers ${\mathbf{N}}_{\bz_p^{*}}(x^3+ax-b)$, ${\mathbf{N}}_{\bz_p}(x^3+ax-b)$, and ${\mathbf{N}}_{\bq_p}(x^3+ax-b)$ of solutions of the depressed cubic equation \eqref{cubiceqn} in domains $\bz_p^{*},$ $\bz_p,$ and $\bq_p,$ respectively, whenever $p>3$ and $ab\neq0.$

Let $a,b\in\bq_p$ be two nonzero $p-$adic numbers with $a=\frac{a^{*}}{|a|_p},$ $b=\frac{b^{*}}{|b|_p}$ where $a^{*},b^{*}\in\bz_p^{*}$ with
$a^{*}=a_0+a_1\cdot p+a_2\cdot p^2+\cdots$ and $b^{*}=b_0+b_1\cdot p+b_2\cdot p^2+\cdots.$

We set $D_0=-4a_0^3-27b_0^2$ and $u_{n+3}=b_0u_n-a_0u_{n+1}$ with $u_1=0,$ $u_2=-a_0,$ and $u_3=b_0$ for $n=\overline{1,p-3}.$

Let $D=-4(a|a|_p)^3-27(b|b|_p)^2$. We have $D=\frac{D^{*}}{|D|_p}$ whenever $D\neq 0,$ where $D^{*}\in\bz_p^{*}$ and
$D^{*}=d_0+d_1\cdot p+d_2\cdot p^2+\cdots$

\begin{thm}\label{numberofsolutionsforp>3} Let $p>3$ be a prime. Then the following statements hold true:
$$
{\mathbf{N}}_{\bz_p^{*}}(x^3+ax-b)=\left\{
\begin{array}{l}
3, \ \ \ |a|_p<|b|_p=1, \ p\equiv 1\ (mod \ 3), \ b_0^{\frac{p-1}{3}}\equiv 1\ (mod \ p)   \\
3, \ \ \ |a|_p=|b|_p=1, \ D=0 \\
3, \ \ \ |a|_p=|b|_p=1, \ 0<|D|_p<1, \ 2\mid\log_p|D|_p, d_0^{\frac{p-1}{2}}\equiv1(mod \ p) \\
3, \ \ \ |a|_p=|b|_p=1, \ |D|_p=1, \ u_{p-2}\equiv 0\ (mod \ p)\\
2, \ \ \ |b|_p<|a|_p=1, \ (-a_0)^{\frac{p-1}{2}}\equiv 1 \ (mod \ p)\\
1, \ \ \ |a|_p<|b|_p=1, \ p\equiv 2\ (mod \ 3) \\
1, \ \ \ |a|_p=|b|_p=1, \ 0<|D|_p<1, \ 2\mid\log_p|D|_p, d_0^{\frac{p-1}{2}}\not\equiv1(mod \ p) \\
1, \ \ \ |a|_p=|b|_p=1, \ 0<|D|_p<1, \ 2\nmid\log_p|D|_p, \\
1, \ \ \ |a|_p=|b|_p=1, \ D_0u_{p-2}^2\not\equiv 0,9a_0^{2}\ (mod \ p)\\
1, \ \ \ |a|_p=|b|_p>1\\
0, \ \ \ otherwise
\end{array}
\right.
$$
$$
{\mathbf{N}}_{\bz_p}(x^3+ax-b)=\left\{
\begin{array}{l}
3, \ \ \ |a|_p^3<|b|_p^2\leq1, \ 3\mid\log_p|b|_p, \ p\equiv 1(mod 3), \ b_0^{\frac{p-1}{3}}\equiv 1\ (mod \ p)   \\
3, \ \ \ |a|_p^3=|b|_p^2\leq1, \ D=0 \\
3, \ \ \ |a|_p^3=|b|_p^2\leq1, \ 0<|D|_p<1, \ 2\mid\log_p|D|_p, d_0^{\frac{p-1}{2}}\equiv1(mod \ p) \\
3, \ \ \ |a|_p^3=|b|_p^2\leq1, \ |D|_p=1, \ u_{p-2}\equiv 0\ (mod \ p)\\
3, \ \ \ |b|_p^2<|a|_p^3\leq1, \ 2\mid\log_p|a|_p, \ (-a_0)^{\frac{p-1}{2}}\equiv 1 \ (mod \ p)\\
1, \ \ \ |a|_p^3<|b|_p^2\leq1, \ 3\mid\log_p|b|_p, \ p\equiv 2\ (mod \ 3) \\
1, \ \ \ |a|_p^3=|b|_p^2\leq1, \ 0<|D|_p<1, \ 2\mid\log_p|D|_p, d_0^{\frac{p-1}{2}}\not\equiv1(mod \ p) \\
1, \ \ \ |a|_p^3=|b|_p^2\leq1, \ 0<|D|_p<1, \ 2\nmid\log_p|D|_p, \\
1, \ \ \ |a|_p^3=|b|_p^2\leq1, \ D_0u_{p-2}^2\not\equiv 0,9a_0^{2}\ (mod \ p)\\
1, \ \ \ |b|_p^2<|a|_p^3\leq1, \ 2\mid\log_p|a|_p, \ (-a_0)^{\frac{p-1}{2}}\not\equiv 1 \ (mod \ p)\\
1, \ \ \ |b|_p^2<|a|_p^3\leq1, \ 2\nmid\log_p|a|_p\\
1, \ \ \ |b|_p^2<|a|_p^3, \ |b|_p\leq|a|_p, \  |a|_p>1\\
0, \ \ \ otherwise
\end{array}
\right.
$$
$$
{\mathbf{N}}_{\bq_p}(x^3+ax-b)=\left\{
\begin{array}{l}
3, \ \ \ |a|_p^3<|b|_p^2, \ 3\mid\log_p|b|_p, \ p\equiv 1(mod 3), \ b_0^{\frac{p-1}{3}}\equiv 1\ (mod \ p)   \\
3, \ \ \ |a|_p^3=|b|_p^2, \ D=0 \\
3, \ \ \ |a|_p^3=|b|_p^2, \ 0<|D|_p<1, \ 2\mid\log_p|D|_p, d_0^{\frac{p-1}{2}}\equiv1(mod \ p) \\
3, \ \ \ |a|_p^3=|b|_p^2, \ |D|_p=1, \ u_{p-2}\equiv 0\ (mod \ p)\\
3, \ \ \ |a|_p^3>|b|_p^2, \ 2\mid\log_p|a|_p, \ (-a_0)^{\frac{p-1}{2}}\equiv 1 \ (mod \ p)\\
1, \ \ \ |a|_p^3<|b|_p^2, \ 3\mid\log_p|b|_p, \ p\equiv 2\ (mod \ 3) \\
1, \ \ \ |a|_p^3=|b|_p^2, \ 0<|D|_p<1, \ 2\mid\log_p|D|_p, d_0^{\frac{p-1}{2}}\not\equiv1(mod \ p) \\
1, \ \ \ |a|_p^3=|b|_p^2, \ 0<|D|_p<1, \ 2\nmid\log_p|D|_p, \\
1, \ \ \ |a|_p^3=|b|_p^2, \ D_0u_{p-2}^2\not\equiv 0,9a_0^{2}\ (mod \ p)\\
1, \ \ \ |a|_p^3>|b|_p^2, \ 2\mid\log_p|a|_p, \ (-a_0)^{\frac{p-1}{2}}\not\equiv 1 \ (mod \ p)\\
1, \ \ \ |a|_p^3>|b|_p^2, \ 2\nmid\log_p|a|_p\\
0, \ \ \ otherwise
\end{array}
\right.
$$
\end{thm}

\begin{pf}
Let $p>3$ be a prime number and $a,b\in\bq_p$ be two nonzero $p-$adic numbers.

{\textsc{Case $\bz_p^{*}$}:} We want to describe the number ${\mathbf{N}}_{\bz_p^{*}}(x^3+ax-b)$ of solutions of the depressed cubic equation \eqref{cubiceqn} in the domain $\bz_p^{*}.$

Due to the case I of Theorem \ref{criterionforp>3}, the depressed cubic equation \eqref{cubiceqn} has a solution in $\bz_p^{*}$ if and only if one of conditions I.1-I.4 holds true. We want to find the number of solutions in every case.

Let us consider the case I.1, i.e., $|a|_p<|b|_p=1$ and $b_0^{\frac{p-1}{(3,p-1)}}\equiv 1 \ (mod \ p).$ In this case, due to Hensel's Lemma \ref{Hensel} as we showed in the proof of I.1, the number of solutions of the equation \eqref{cubiceqn} in $\bz_p^{*}$ is the same as the number of solutions of the equation $x_0^3= b_0$ in $\bbf_p.$ Then, due to Proposition \ref{aisresidueofp}, the last equation has 3 distinct solutions if $p\equiv 1 \ (mod \ 3)$ and it has a unique solution if $p\equiv 2 \ (mod \ 3).$

Therefore, if $|a|_p<|b|_p=1$, $p\equiv 1 \ (mod \ 3)$, and $b_0^{\frac{p-1}{3}}\equiv 1 \ (mod \ p)$ then the equation \eqref{cubiceqn} has 3 distinct solutions in $\bz_p^{*}$ and if $|a|_p<|b|_p=1,$ $p\equiv 2 \ (mod \ 3)$ and $b_0^{p-1}\equiv 1 \ (mod \ p)$ then the equation \eqref{cubiceqn} has a unique solution $\bz_p^{*}.$ It is worth mentioning that one always has $b_0^{p-1}\equiv 1 \ (mod \ p)$ since $(b_0,p)=1$.

Let us consider the case I.2, i.e., $|b|_p<|a|_p=1$ and $(-a_0)^{\frac{p-1}{2}}\equiv 1 \ (mod \ p).$ In this case, due to Hensel's Lemma \ref{Hensel} as we showed in the proof of I.2, the number of solutions of the equation \eqref{cubiceqn} in $\bz_p^{*}$ is the same as the number of solutions of the equation $x_0^2=-a_0$ in $\bbf_p.$ Since $(-a_0)^{\frac{p-1}{2}}\equiv 1 \ (mod \ p),$ the last equation has 2 distinct solutions in $\bbf_p.$

Therefore, if $|b|_p<|a|_p=1$ and $(-a_0)^{\frac{p-1}{2}}\equiv 1 \ (mod \ p)$ then the depressed cubic equation \eqref{cubiceqn} has 2 distinct solutions in $\bz_p^{*}$.

Let us consider the case I.3, i.e., $|a|_p=|b|_p=1$ and $D_0u_{p-2}^2\not\equiv 9a_0^{2} \ (mod \ p).$ In this case, as we showed in the proof of I.3, the equation \eqref{cubiceqn} has a solution $\bar{x}$ such that $\bar{x}\equiv\bar{x}_0 \ (mod \ p),$ where $\bar{x}_0$ is a solution of the following congruent equation
\begin{eqnarray}\label{I3}
x_0^3+a_0x_0\equiv b_0 \ (mod \ p)
\end{eqnarray}
such that $3\bar{x}_0^2+a_0\not\equiv 0 \ (mod \ p).$

Due to Proposition \ref{CubicinF_p}, if $D_0u_{p-2}^2\not\equiv 0, 9a_0^{2} \ (mod \ p)$ then the equation \eqref{I3} does not have any solution except $\bar{x}_0$. Therefore, due to Hensel's Lemma \ref{Hensel}, if $|a|_p=|b|_p=1$ and $D_0u_{p-2}^2\not\equiv 0, 9a_0^{2} \ (mod \ p)$ then the equation \eqref{cubiceqn} has a unique solution in $\bz_p^{*}$.

Now, let us study the case $|a|_p=|b|_p=1$ and $D_0u_{p-2}^2\equiv 0 \ (mod \ p)$. In this case, due to Proposition \ref{CubicinF_p}, the congruent equation \eqref{I3} has 2 more solutions besides $\bar{x}_0.$ We denote them by $x_0$ and $y_0$. There is no loss of generality in assuming that $0<x_0,y_0<p.$

Due to Proposition \ref{numberofcongequation}, if $D_0\not\equiv 0 \ (mod \ p)$ then all solutions $x_0,y_0,\bar{x}_0$ of the congruent equation \eqref{I3} are distinct from each other and $3{x}_0^2+a_0\not\equiv 0 \ (mod \ p),$ $3{y}_0^2+a_0\not\equiv 0 \ (mod \ p).$ Therefore, due to Hensel's Lemma \ref{Hensel}, the equation \eqref{cubiceqn} has 2 more solutions $x,$ $y$ besides $\bar{x}$ such that ${x}\equiv{x}_0 \ (mod \ p),$ ${y}\equiv{y}_0 \ (mod \ p).$ Hence, $|a|_p=|b|_p=1,$  $D_0\not\equiv 0 \ (mod \ p)$ (equivalently $|D|_p=1$) and $u_{p-2}\equiv 0 \ (mod \ p)$ then the equation \eqref{cubiceqn} has 3 distinct solutions in $\bz_p^{*}$.

Again according to Proposition \ref{numberofcongequation}, if $D_0\equiv 0 \ (mod \ p)$ then two solutions $x_0,y_0$ of the congruent equation \eqref{I3} are equal to each other and $\bar{x}_0=-2x_0,$ $3{x}_0^2+a_0\equiv 0 \ (mod \ p).$ Now, we are going to study this case in a detail.

Since $\bar{x}$ is a solution of the equation \eqref{cubiceqn}, one can get that
$$
x^3+ax-b=(x-\bar{x})(x^2+\bar{x}x+\bar{x}^2+a).
$$

We are going to study the following quadratic equation
\begin{eqnarray}\label{quadraticequat}
x^2+\bar{x}x+\bar{x}^2+a=0.
\end{eqnarray}
It is clear that
$$
\left(x+\frac{\bar{x}}{2}\right)^2=-\left(3\left(\frac{\bar{x}}{2}\right)^2+a\right).
$$
Since $a,\bar{x}\in\bz_p^{*}$ and $p>3$, we have that $\left|3\left(\frac{\bar{x}}{2}\right)^2+a\right|_p\leq1.$ Since $\bar{x}$ is a solution of the equation \eqref{cubiceqn}, one can get that
$
4\bar{x}\left(3\left(\frac{\bar{x}}{2}\right)^2+a\right)=a\bar{x}+3b.
$
We then  have that
\begin{eqnarray*}
D&=&-4a^3-27b^2=3(a^2\bar{x}^2-9b^2)-a^2(3\bar{x}^2+4a)\\
&=&12\bar{x}\left(3\left(\frac{\bar{x}}{2}\right)^2+a\right)(a\bar{x}-3b)-4a^2\left(3\left(\frac{\bar{x}}{2}\right)^2+a\right)\\
&=&4\left(3\left(\frac{\bar{x}}{2}\right)^2+a\right)\left(3a\bar{x}^2-9b\bar{x}-a^2\right).
\end{eqnarray*}

Since $\bar{x}\equiv\bar{x}_0=-2x_0 \ (mod \ p),$ $3x_0^2+a_0\equiv 0 \ (mod \ p)$, and $2x_0a_0\equiv 3b_0\ (mod \ p)$ we get that $3a\bar{x}^2-9b\bar{x}-a^2\equiv -9a_0^2\ (mod \ p).$ This means that $\left|3a\bar{x}^2-9b\bar{x}-a^2\right|_p=1$ and $-\left(3a\bar{x}^2-9b\bar{x}-a^2\right)$ is a complete square of some $p-$adic integer number. Therefore, we obtain that
\begin{eqnarray}\label{Dand3x_0^2+a}
-\left(3\left(\frac{\bar{x}}{2}\right)^2+a\right)=\frac{D}{-4\left(3a\bar{x}^2-9b\bar{x}-a^2\right)}
\end{eqnarray}

Let us analyze the quadratic equation \eqref{quadraticequat}.

If $D=0$ then the quadratic equation \eqref{quadraticequat} has solutions $x_1=x_2=\frac{3b}{2a}$ in $\bz_p^{*}.$ Therefore, $|a|_p=|b|_p=1$ and $D=0$ then the cubic equation \eqref{cubiceqn} has 3 solutions such that $x_1=x_2=\frac{3b}{2a}$ and $x_3=-\frac{3b}{a.}$

Let $D\neq 0.$ Since $D\equiv D_0\equiv 0 \ (mod \ p)$, there exists $k\in\bn$ such that $|D|_p=p^{-k}$, i.e., $D=\frac{D^{*}}{|D|_p},$ where $D^{*}\in\bz_p^{*}$ with $D^{*}=d_0+d_1\cdot p+d_2\cdot p^2+\cdots$.

The quadratic equation \eqref{quadraticequat} has a solution if and only if $\log_p|D|_p=-k$ is even number and $d_0^{\frac{p-1}{2}}\equiv 1 \ (mod \ p)$. In this case $-\left(3\left(\frac{\bar{x}}{2}\right)^2+a\right)$ is a complete square and the quadratic equation \eqref{quadraticequat} has  two distinct solutions in $\bz_p^{*}$ as follows
\begin{eqnarray}
x_{\pm}=-\frac{\bar{x}}{2}\pm\sqrt{-\left(3\left(\frac{\bar{x}}{2}\right)^2+a\right)}
\end{eqnarray}

Therefore, if $|a|_p=|b|_p=1,$ $0<|D|_p<1,$ $2\mid \log_p|D|_p,$ and $d_0^{\frac{p-1}{2}}\equiv 1 \ (mod \ p)$ then the cubic equation \eqref{cubiceqn} has 3 distinct solutions in $\bz_p^{*}.$ If $|a|_p=|b|_p=1,$ $0<|D|_p<1,$ $2\mid \log_p|D|_p,$ and $d_0^{\frac{p-1}{2}}\not\equiv 1 \ (mod \ p)$ or $|a|_p=|b|_p=1,$ $0<|D|_p<1,$ and $2\nmid \log_p|D|_p$ then the cubic equation \eqref{cubiceqn} has a unique solutions in $\bz_p^{*}.$

Let us consider the case I.4, i.e., $|a|_p=|b|_p>1.$ In this case, due to Hensel's Lemma \ref{Hensel} as we showed in the proof of I.4, the number of solutions of the equation \eqref{cubiceqn} in $\bz_p^{*}$ is the same as the number of solutions of the linear equation $a_0x_0=b_0$ in $\bbf_p.$ Since $a_0\neq 0$, the last equation has a unique solution. Therefore, if  $|a|_p=|b|_p>1$ then the depressed cubic equation \eqref{cubiceqn} has a unique solution in $\bz_p^{*}$.

{\textsc{Case $\bz_p$}:} We shall study the number ${\mathbf{N}}_{\bz_p}(x^3+ax-b)$ of solutions of the cubic equation \eqref{cubiceqn} in the domain $\bz_p.$

Due to the case II of Theorem \ref{criterionforp>3}, the cubic equation \eqref{cubiceqn} has a solution in $\bz_p$ if and only if one of conditions II.1-II.3 holds true. We want to find the number of solutions in every case.

Let us consider the case II.1, i.e., $|a|_p^3<|b|_p^2\leq 1,$  $3\mid\log_p|b|_p,$ and $b_0^{\frac{p-1}{(3,p-1)}}\equiv 1 \ (mod \ p).$ In this case, as we showed in the proof of II.1, the number of  solutions of the cubic equation \eqref{cubiceqn} in the domain $\bz_p$ is the same as the number of solutions of the following equation in the domain $\bz_p^{*}:$
\begin{eqnarray}\label{II.1}
y^3+a\sqrt[3]{|b|_p^2}y=b^{*}.
\end{eqnarray}

Then it is clear that $\left|a\sqrt[3]{|b|_p^2}\right|_p<\left|b^{*}\right|_p=1$ and $b_0^{\frac{p-1}{(3,p-1)}}\equiv 1 \ (mod \ p).$ In this case, as we already discussed that if $p\equiv 1 \ (mod \ 3)$ then the equation \eqref{II.1} has 3 distinct solutions in $\bz_p^{*}$ and if $p\equiv 2 \ (mod \ 3)$ then the equation \eqref{II.1} has a unique solution in $\bz_p^{*}$.

Consequently, if $|a|_p^3<|b|_p^2\leq 1,$  $3\mid\log_p|b|_p,$ $p\equiv 1 \ (mod \ 3)$, and $b_0^{\frac{p-1}{3}}\equiv 1 \ (mod \ p)$ then the cubic equation \eqref{cubiceqn} has 3 distinct solutions in $\bz_p$ and if $|a|_p^3<|b|_p^2\leq 1,$  $3\mid\log_p|b|_p,$ and $p\equiv 2 \ (mod \ 3)$ then the cubic equation \eqref{cubiceqn} has a unique solution in $\bz_p.$

Let us consider the case II.2, i.e., $|a|_p^3=|b|_p^2\leq 1$ and $D_0u_{p-2}^2\not\equiv 9a_0^{2} \ (mod \ p).$ In this case, as we showed in the proof of II.2 that the number of  solutions of the cubic equation \eqref{cubiceqn} in the domain $\bz_p$ is the same as the number of solutions of the following equation in the domain $\bz_p^{*}$
\begin{eqnarray}\label{II.2}
y^3+a^{*}y=b^{*}.
\end{eqnarray}

Then it is clear that $|a^{*}|_p=|b^{*}|_p=1$ and $D_0u_{p-2}^2\not\equiv 9a_0^{2} \ (mod \ p).$

Let $D=-4(a^{*})^3-27(b^{*})^2$. We have $D=\frac{D^{*}}{|D|_p}$ whenever $D\neq 0,$ where $D^{*}\in\bz_p^{*}$ and
$D^{*}=d_0+d_1\cdot p+d_2\cdot p^2+\cdots.$

The number of solutions of the equation \eqref{II.2} was studied very well in the case $|a^{*}|_p=|b^{*}|_p=1$ and $D_0u_{p-2}^2\not\equiv 9a_0^{2} \ (mod \ p).$ The equation \eqref{II.2} has either a unique solution or 3 solutions.

Consequently, the cubic equation \eqref{cubiceqn} has 3 solutions if and only if one of the following conditions holds true: (i) $|a|_p^3=|b|_p^2\leq 1$ and $D=0$ or (ii) $|a|_p^3=|b|_p^2\leq 1,$ $0<|D|_p<1,$ $2\mid\log_p|D|_p,$ and $d_0^{\frac{p-1}{2}}\equiv1 \ (mod \ p)$ or (iii) $|a|_p^3=|b|_p^2\leq 1,$ $|D|_p=1,$ and $u_{p-2}\equiv 0 \ (mod \ p)$. The cubic equation \eqref{cubiceqn} has a unique solution if and only if one of the following conditions holds true: (i) $|a|_p^3=|b|_p^2\leq 1,$ $0<|D|_p<1,$ $2\mid\log_p|D|_p,$ and $d_0^{\frac{p-1}{2}}\not\equiv1 \ (mod \ p)$ or (ii) $|a|_p^3=|b|_p^2\leq 1,$ $0<|D|_p<1,$ and $2\nmid\log_p|D|_p$ or (iii) $|a|_p^3=|b|_p^2\leq 1$ and $D_0u_{p-2}^2\not\equiv 0,9a_0^{2} \ (mod \ p).$

Let us consider the case II.3, i.e., $|a|_p^3>|b|_p^2$ and $|a|_p\geq |b|_p.$

One can easily check that
$$
\Delta=\Delta_1\cup\Delta_2\cup\Delta_3\cup\Delta_4,\\
$$
where
\begin{eqnarray*}
\Delta&=&\left\{(a,b): |a|_p^3>|b|_p^2, \ |a|_p\geq |b|_p\right\},\\
\Delta_1&=&\left\{(a,b): |b|_p^2<|a|_p^3\leq 1, \ 2\mid\log_p|a|_p, (-a_0)^{\frac{p-1}{2}}\equiv 1 \ (mod \ p)\right\},\\
\Delta_2&=&\left\{(a,b): |b|_p^2<|a|_p^3\leq 1, \ 2\mid\log_p|a|_p, (-a_0)^{\frac{p-1}{2}}\not\equiv 1 \ (mod \ p)\right\},\\
\Delta_3&=&\left\{(a,b): |b|_p^2<|a|_p^3\leq 1, \ 2\nmid\log_p|a|_p \right\},\\
\Delta_4&=&\{(a,b): |a|_p^3>|b|_p^2, \ |a|_p\geq |b|_p,\ |a|_p>1\}.
\end{eqnarray*}

In the case $\Delta_1$, as we showed in the proof of II.3, the number of  solutions of the cubic equation \eqref{cubiceqn} in the domain $\bz_p$ is the same as the total number of solutions of the following equations in the domain $\bz_p^{*}:$
\begin{eqnarray}
\label{II.3.1}y^3+a^{*}y=b\sqrt{|a|_p^3}, \\
\label{II.3.2}z^3+a\left|\frac{b}{a}\right|_p^2z=b\left|\frac{b}{a}\right|_p^3,
\end{eqnarray}
It is then clear that $|a^{*}|_p=1,$ $\left|b\sqrt{|a|_p^3}\right|_p<1$ and $\left|a\left|\frac{b}{a}\right|_p^2\right|_p=\left|b\left|\frac{b}{a}\right|_p^3\right|_p>1$.

In this case, as we already discussed, the equation \eqref{II.3.1} has 2 distinct solutions in $\bz_p^{*}$ and the equation \eqref{II.3.2} has a unique solution in $\bz_p^{*}$. Consequently, the cubic equation \eqref{cubiceqn} has 3 solutions in $\bz_p$.

In the case $\Delta_2\cup\Delta_3\cup\Delta_4$, as we showed in the proof of II.3, the number of  solutions of the cubic equation \eqref{cubiceqn} in the domain $\bz_p$ is the same as the number of solutions of the equation \eqref{II.3.2} in the domain $\bz_p^{*}.$ We know that the equation \eqref{II.3.2} has a unique solution in $\bz_p^{*}$. Consequently, the  cubic equation \eqref{cubiceqn} has a unique solution in $\bz_p$.

{\textsc{Case $\bq_p$}:} Analogously, one can study the number ${\mathbf{N}}_{\bq_p}(x^3+ax-b)$ of solutions of the cubic equation \eqref{cubiceqn} in the domain $\bq_p$. This completes the proof.

\end{pf}

%

\section{A $p-$adic Cardano formula}

As we already mentioned, in general, by means of Cardano's method we could not detect solutions of all cubic equations over the $p-$adic filed. In this section, we shall describe all possible cases in which Cardano's method is applicable to solve the cubic equation in $\bq_p$.

Let $a,b\in\bq_p$ be two nonzero $p-$adic numbers with $a=\frac{a^{*}}{|a|_p},$ $b=\frac{b^{*}}{|b|_p}$ where $a^{*},b^{*}\in\bz_p^{*}$ with
$a^{*}=a_0+a_1\cdot p+a_2\cdot p^2+\cdots$ and $b^{*}=b_0+b_1\cdot p+b_2\cdot p^2+\cdots.$

We set $D_0=-4a_0^3-27b_0^2$ and $u_{n+3}=b_0u_n-a_0u_{n+1}$ with $u_1=0,$ $u_2=-a_0,$ and $u_3=b_0$ for $n=\overline{1,p-3}.$

Let $D=-4(a|a|_p)^3-27(b|b|_p)^2$. We have $D=\frac{D^{*}}{|D|_p}$ whenever $D\neq 0,$ where $D^{*}\in\bz_p^{*}$ and
$D^{*}=d_0+d_1\cdot p+d_2\cdot p^2+\cdots$

Let $\Delta_0\in\{1,2,\cdots,p-1\}$ such that $\Delta_0^2\equiv -3d_0\  (mod\ p)$ whenever $-3d_0$ is a quadratic residue.

\begin{thm}
The Cardano method is applicable for the depressed cubic equation \eqref{cubiceqn} in $\bq_p$
and one of solutions has the following form
\begin{eqnarray}\label{Cardano form}
x=\sqrt[3]{\frac{b}{2}+\sqrt{\left(\frac{a}{3}\right)^3+\left(\frac{b}{2}\right)^2}}+
\sqrt[3]{\frac{b}{2}-\sqrt{\left(\frac{a}{3}\right)^3+\left(\frac{b}{2}\right)^2}}
\end{eqnarray}
if and only if one of the following conditions holds true:
\begin{itemize}
    \item [I.1] $|a|_p^3<|b|_p^2,$ $3\mid\log_p|b|_p,$ and $b_0^{\frac{p-1}{(3,p-1)}}\equiv 1 \ (mod \ p);$
    \item [I.2] $|a|_p^3=|b|_p^2,$ $D=0,$ and $(4b_0)^{\frac{p-1}{(3,p-1)}}\equiv 1 \ (mod \ p);$
    \item [I.3] $|a|_p^3=|b|_p^2,$ $0<|D|_p<1,$ $2\mid\log_p|D|_p,$ \\
    $(-3d_0)^{\frac{p-1}{2}}\equiv 1 \ (mod \ p ),$ and $(4b_0)^{\frac{p-1}{(3,p-1)}}\equiv 1 \ (mod \ p);$
    \item [I.4] $|a|_p^3=|b|_p^2,$ $|D|_p=1,$ $D_0u_{p-2}^2\not\equiv 9a_0^{2} \ (mod \ p),$\\
    $(-3d_0)^{\frac{p-1}{2}}\equiv 1 \ (mod \ p ),$ and $(108b_0+12\Delta_0)^{\frac{p-1}{(3,p-1)}}\equiv 1 \ (mod \ p);$
    \item [I.5] $|a|_p^3>|b|_p^2,$ $2\mid\log_p|a|_p,$ and $(3a_0)^{\frac{p-1}{2}}\equiv 1 \ (mod \ p )$
  \end{itemize}
\end{thm}

\begin{pf}
We want to describe all $a,b\in\bq_p$ in which the expression \eqref{Cardano form} is well defined in $\bq_p.$
Let us consider III case of Theorem \ref{criterionforp>3} where the cubic equation has at least one solution.

III.1  $|a|_p^3<|b|_p^2,$  $3\mid\log_p|b|_p,$ and $b_0^{\frac{p-1}{(3,p-1)}}\equiv 1 \ (mod \ p).$

Now, we show that the expression \eqref{Cardano form} is meaningful in $\bq_p.$ In fact, since $\left|\left(\frac{a}{3}\right)^3\right|_p<\left|\left(\frac{b}{2}\right)^2\right|_p$, the expression $\sqrt{\left(\frac{a}{3}\right)^3+\left(\frac{b}{2}\right)^2}$ is well defined and $\left|\sqrt{\left(\frac{a}{3}\right)^3+\left(\frac{b}{2}\right)^2}\right|_p=\left|\frac{b}{2}\right|_p=|b|_p.$

Since $3\mid\log_p|b|_p,$ and $b_0^{\frac{p-1}{(3,p-1)}}\equiv 1 \ (mod \ p),$ the expression   $\frac{b}{2}+\sqrt{\left(\frac{a}{3}\right)^3+\left(\frac{b}{2}\right)^2}$ is a cube of some $p-$adic number. It means that $\sqrt[3]{\frac{b}{2}+\sqrt{\left(\frac{a}{3}\right)^3+\left(\frac{b}{2}\right)^2}}$ is meaningful. On the other hand, we have that
\begin{eqnarray}\label{multipexpression}
\left(\frac{b}{2}+\sqrt{\left(\frac{a}{3}\right)^3+\left(\frac{b}{2}\right)^2}\right)
\left(\frac{b}{2}-\sqrt{\left(\frac{a}{3}\right)^3+\left(\frac{b}{2}\right)^2}\right)=-\left(\frac{a}{3}\right)^3.
\end{eqnarray}

Therefore, the expressions  $\sqrt[3]{\frac{b}{2}-\sqrt{\left(\frac{a}{3}\right)^3+\left(\frac{b}{2}\right)^2}}$ and \eqref{Cardano form} are meaningful.

III.2 $|a|_p^3=|b|_p^2$ and $D_0u_{p-2}^2\not\equiv 9a_0^{2} \ (mod \ p);$

It is clear that $\left(\frac{a}{3}\right)^3+\left(\frac{b}{2}\right)^2=\frac{1}{324|b|_p^2}\cdot (-3D).$ Then the last expression is a perfect square if and only if one of the following conditions holds true:
\begin{itemize}
  \item [(i)] $D=0;$
  \item [(ii)] $0<|D|_p<1,$ $2\mid\log_p|D|_p,$ and $(-3d_0)^{\frac{p-1}{2}}\equiv 1 \ (mod \ p );$
  \item [(iii)] $|D|_p=1$ and $(-3d_0)^{\frac{p-1}{2}}\equiv 1 \ (mod \ p ).$
\end{itemize}

We shall study every case.

Let $D=0.$ Then the expression \eqref{Cardano form} is meaningful if and only if $\frac{b}{2}$ is a cube of some $p-$adic number, i.e., $(4b_0)^{\frac{p-1}{(3,p-1)}}\equiv 1 \ (mod \ p).$

Let $0<|D|_p<1,$ $2\mid\log_p|D|_p,$ and $(-3d_0)^{\frac{p-1}{2}}\equiv 1 \ (mod \ p ).$

In this case, we have that $\left|\frac{b}{2}\right|_p>\left|\sqrt{\left(\frac{a}{3}\right)^3+\left(\frac{b}{2}\right)^2}\right|_p.$ Therefore, the expression \eqref{Cardano form} is meaningful if and only if $\frac{b}{2}$ is a cube of some $p-$adic number, i.e., $(4b_0)^{\frac{p-1}{(3,p-1)}}\equiv 1 \ (mod \ p).$

Let $|D|_p=1$ and $(-3d_0)^{\frac{p-1}{2}}\equiv 1 \ (mod \ p ).$

In this case, we obtain that $\left|\frac{b}{2}\right|_p=\left|\sqrt{\left(\frac{a}{3}\right)^3+\left(\frac{b}{2}\right)^2}\right|_p.$ Due to \eqref{multipexpression}, we get that
$$
\left|\frac{b}{2}\pm\sqrt{\left(\frac{a}{3}\right)^3+\left(\frac{b}{2}\right)^2}\right|=\left|\frac{b}{2}\right|_p.
$$

Therefore, the expression \eqref{Cardano form} is meaningful if and only if $(108b_0+12\Delta_0)^{\frac{p-1}{(3,p-1)}}\equiv 1 \ (mod \ p),$ where $\Delta_0\in\{1,2,\cdots,p-1\}$ such that $\Delta_0^2\equiv -3d_0\  (mod\ p).$

III.3 $|a|_p^3>|b|_p^2.$

We know that $\left|\left(\frac{a}{3}\right)^3\right|_p>\left|\left(\frac{b}{2}\right)^2\right|_p$. Then, the expression $\sqrt{\left(\frac{a}{3}\right)^3+\left(\frac{b}{2}\right)^2}$ is well defined if and only if $\frac{a}{3}$ is a perfect square. The last one means that $2\mid\log_p|a|_p,$ and $(3a_0)^{\frac{p-1}{2}}\equiv 1 \ (mod \ p ).$

In this case, we have that $\left|\left(\sqrt{\frac{a}{3}}\right)^3\right|_p>\left|\frac{b}{2}\right|_p.$ Therefore, the expressions ${\frac{b}{2}+\sqrt{\left(\frac{a}{3}\right)^3+\left(\frac{b}{2}\right)^2}}$ and ${\frac{b}{2}-\sqrt{\left(\frac{a}{3}\right)^3+\left(\frac{b}{2}\right)^2}}$ are cube of some $p-$adic numbers. This means that the expression \eqref{Cardano form} is well defined in $\bq_p.$

\end{pf}

\section*{Acknowledgement}  The Author (M.S.) is grateful to Pah Chin Hee for his valued discussion. We are thank Jan Kohlhaase for his attention and some comments on our papers.

\end{document}